\documentclass[12pt]{article}
\usepackage{graphicx,xspace,colortbl}
\usepackage{amsmath,amsthm,amsfonts}
\usepackage{authblk}
\usepackage{color}
\usepackage{fancybox}
\usepackage{amssymb,mathrsfs}
\usepackage{pdfsync}
\RequirePackage[colorlinks,citecolor=blue,urlcolor=blue,linkcolor=blue]{hyperref}
\hypersetup{
colorlinks = true,
citecolor=blue,
urlcolor=blue,
linkcolor=blue,
pdfauthor = {Peng-Cheng Hang, Alexey Kuznetsov},
pdfkeywords = {},
pdftitle = {Fourier transform pairs and Eisenstein-type series related to Jacobi elliptic functions},
pdfpagemode = UseNone
}
    \def\qed{\hfill$\sqcap\kern-8.0pt\hbox{$\sqcup$}$\\}

    \def\numx#1e#2{{#1}\mathrm{e}{#2}}
    \def\re{\textnormal {Re}}
    \def\im{\textnormal {Im}}
    \def\r{{\mathbb R}}
    \def\c{{\mathbb C}}
    \def\d{{\textnormal d}}
    \def\i{{\textnormal i}}
    \def\mi{{\textnormal i}}

 	\newtheorem{theorem}{Theorem}
	
	\newtheorem{proposition}{Proposition}

\newcommand{\CC}{\mathbb{C}}

\newcommand{\RR}{\mathbb{R}}

\newcommand{\ZZ}{\mathbb{Z}}

\newcommand{\md}{\mathrm{d}}

\newcommand{\mo}{\mathcal{O}}
\newcommand{\sm}{\setminus}

\newcommand{\ns}{\mathrm{ns}}
\newcommand{\nc}{\mathrm{nc}}
\newcommand{\nd}{\mathrm{nd}}
\newcommand{\sn}{\mathrm{sn}}
\newcommand{\sd}{\mathrm{sd}}
\newcommand{\cn}{\mathrm{cn}}
\newcommand{\cs}{\mathrm{cs}}
\newcommand{\cd}{\mathrm{cd}}
\newcommand{\dn}{\mathrm{dn}}
\newcommand{\dc}{\mathrm{dc}}
\newcommand{\ds}{\mathrm{ds}}

\renewcommand{\sc}{\mathrm{sc}}

\DeclareMathOperator*{\Res}{Res}

\title{Fourier transform pairs and Eisenstein-type series related to Jacobi elliptic functions}
\author{Peng-Cheng Hang\thanks{School of Mathematics and Statistics, Donghua University, Shanghai 201620, P.R. China \\ Email: mathroc618@outlook.com}, Alexey Kuznetsov\thanks{Department of Mathematics and Statistics, York University, Toronto, Ontario, M3J 1P3, Canada \\ Email: akuznets@yorku.ca}}

\date{\today}
 
\begin{document}
\maketitle

\begin{abstract}
We compute Fourier transforms of functions expressed as a ratio of one of the Jacobi elliptic functions divided by $\sinh(\pi x)$ or $\cosh(\pi x)$. In many cases, the resulting Fourier transform remains within the same class of functions. Applying the Mellin transform, we obtain sixteen Eisenstein-type series $\zeta_{j,l}(s,\tau)$, for which we establish several results: analytic continuation with respect to the variable $s$, a functional equation connecting $\zeta_{j,l}(s,\tau)$ and $\zeta_{l,j}(1-s,-1/\tau)$, and explicit expressions for $\zeta_{j,l}(s,\tau)$ when $s$ runs through a sequence of positive even or odd integers.
	\vspace{4mm}
	
	\vspace{2mm}
	
\end{abstract}

\section{Introduction}

Consider the following integral \cite[\S 3.981]{Jeffrey2007}:
\begin{equation}\label{Fs_sinh}
\int_0^{\infty} \frac{\sin(xy)}{\sinh(x)}\d x=
\frac{\pi}{2} \tanh\Big(\frac{\pi y}{2}\Big), \;\;\; y \in \r.  
\end{equation}
This result has several important consequences. The Mellin transform of $1/\sinh(x)$ can be found in \cite[\S 3.523]{Jeffrey2007} and is given by
\begin{equation}\label{sinh_integral_zeta}
\int_0^{\infty} \frac{x^{s-1} \d x}{\sinh(x)} =2 (1-2^{-s}) \Gamma(s) \zeta(s), \;\;\; \re(s)>1.
\end{equation}
Taking the derivative $\d^{2n-1}/\d y^{2n-1}$ of both sides of \eqref{Fs_sinh} and setting $y=0$ yields formulas for $\zeta(2n)$, $n \in {\mathbb N}$, expressed in terms of derivatives of $\tanh(x)$, which can in turn be written using the Bernoulli numbers $B_{2n}$; see \cite[\S 1.411]{Jeffrey2007}. The Mellin transform of the function $\tanh(x)$ is also expressed in terms of the Riemann zeta function (see formula (1.6.3) in \cite{Ober1}). Using this result and the identity \eqref{Fs_sinh}, one can apply the technique from \cite[Section 2.7]{TiHB} to derive the functional equation for the Riemann zeta function, which states that the function  
\begin{equation}\label{def_Lambda_zeta}
\Lambda(s):=\pi^{-\frac{s}{2}}\Gamma\Big(\frac{s}{2}\Big)\zeta(s)
\end{equation}
is invariant under the transformation $s\mapsto 1-s$. 

A similarly important integral is 
\begin{equation}\label{Fc_cosh}
\int_0^{\infty} \frac{\cos(xy)}{\cosh(x)}\d x=
 \frac{\pi/2}{\cosh(\pi y/2)}, \;\;\; y \in \r.  
\end{equation} 
The Mellin transform of $1/\cosh(x)$ is
\[
\int_0^{\infty} \frac{x^{s-1}\d x}{\cosh(x)}=2\Gamma(s)\beta(s),\ \ \re(s)>0,
\]
where 
\[
\beta(s)=L(s,\chi_4)=\sum_{n\geqslant 0}\frac{(-1)^n}{(2n+1)^s},\;\;\; \re(s)>0,
\]
is the Dirichlet beta function (equivalently, the Dirichlet $L$-function associated with the primitive Dirichlet character modulo $4$).  
These results can be found in \cite[\S 3.981]{Jeffrey2007} and \cite[\S 3.523]{Jeffrey2007}. Taking the derivative $\d^{2n}/\d y^{2n}$ of both sides of \eqref{Fc_cosh} and setting $y=0$ yields explicit expressions for the values $\beta(2n+1)$, $n\in {\mathbb Z}_{\geqslant 0}$, given in terms of the Euler numbers $E_{2n}$. Moreover, the invariance of $1/\cosh(x)$ under the Fourier cosine transform in \eqref{Fc_cosh} leads to the functional equation for the Dirichlet beta function: the function
\begin{equation}\label{def_Lambda_beta}
\Lambda(s,\chi_4):=2^s\pi^{-\frac{s}{2}}\Gamma\Big(\frac{s+1}{2}\Big)\beta(s)
\end{equation}
is invariant under the transformation $s\mapsto 1-s$.

In this paper, we pursue two main objectives. First, we present additional examples of meromorphic periodic functions with explicit Fourier transforms. Second, by applying the Mellin transform to these functions, we obtain a number of double series possessing interesting analytic properties.

\begin{table}[t]
	\centering
    \vspace{1mm}
	\begin{tabular}{|c|c|c||c|c|c|}
		\hline
		 & $f(x)$ & ${\mathcal F}[f](y)$ & & $f(x)$ & ${\mathcal F}[f](y)$   \\
		\hline
		\rule[-7.5pt]{0pt}{32pt}
		1& $\displaystyle \frac{\dc(2K' x,k )}{\cosh(\pi x)}$ &  $\displaystyle -2k\,\frac{k'\, \sd(2Ky,k')-e^{\pi y}}{e^{2\pi y}+1}$ & 10 & 
		$\displaystyle \frac{\dc(2K' x,k)}{\sinh(\pi x)}$ & $\displaystyle   2\mi\,\frac{k' \, \cd(2Ky,k')-1}{e^{2\pi y}+1}$ \\
		\rule[-7.5pt]{0pt}{32pt}
		2 & $\displaystyle \frac{\nc(2K' x,k )}{\cosh(\pi x)}$ & $\displaystyle 2k\,\frac{\sd(2Ky,k')}{e^{2\pi y}-1}$ & 11 &
		$\displaystyle \frac{\nc(2K' x,k )}{\sinh(\pi x)}$ & $\displaystyle  -2\mi\,\frac{\cd(2Ky,k')-1}{e^{2\pi y}-1}$ \\
		\rule[-7.5pt]{0pt}{32pt} 
		3 & $\displaystyle \frac{\sc(2K' x,k)}{\cosh(\pi x)}$ & $\displaystyle -2\mi\,\frac{\cn(2Ky,k')-e^{\pi y}}{e^{2\pi y}-1}$ & 12 &
		$\displaystyle \frac{\sc(2K' x,k)}{\sinh(\pi x)}$ & $\displaystyle  2\,\frac{\sn(2Ky,k')}{e^{2\pi y}-1}$ \\ 
		\rule[-7.5pt]{0pt}{32pt}
		4 & $\displaystyle \frac{\cd(2K' x,k )}{\cosh(\pi x)}$ &  $\displaystyle -\frac{1}{k}\,\frac{k' \; \cd(2Ky,k')-1}{\cosh(\pi y)}$ & 13 &
		$\displaystyle \frac{\cd(2K' x,k )}{\sinh(\pi x)}$ & $\displaystyle  -\mi\,\frac{k' \, \sd(2Ky,k')+e^{-\pi y}}{\cosh(\pi y)}$ \\
		\rule[-7.5pt]{0pt}{32pt} 
		5 & $\displaystyle \frac{\nd(2K' x,k )}{\cosh(\pi x)}$ &  $\displaystyle \frac{\sn(2Ky,k')}{\sinh(\pi y)}$ & 14 &
		$\displaystyle \frac{\nd(2K' x,k )}{\sinh(\pi x)}$ &  $\displaystyle -\mi\,\frac{\cn(2Ky,k')-e^{-\pi y}}{\sinh(\pi y)}$ \\ 
		\rule[-7.5pt]{0pt}{32pt} 
		6 & $\displaystyle \frac{\sd(2K' x,k )}{\cosh(\pi x)}$ & $\displaystyle  -\frac{\mi}{k}\,\frac{\cd(2Ky,k')-1}{\sinh(\pi y)}$
		& 15 &
		$\displaystyle \frac{\sd(2K' x,k )}{\sinh(\pi x)}$ & $\displaystyle  \frac{\sd(2Ky,k')}{\sinh(\pi y)}$ \\
		\rule[-7.5pt]{0pt}{32pt}
		7 & $\displaystyle \frac{\cs(2K' x,k )}{\cosh(\pi x)}$ &  $\displaystyle  2\mi\,\frac{\dn(2Ky,k')-e^{\pi y}}{e^{2\pi y}-1}$ 
		& 16 &
		$\displaystyle \frac{\cn(2K' x,k)}{\sinh(\pi x)}$ & $\displaystyle  -\mi\,\frac{\nd(2Ky,k')-e^{-\pi y}}{\sinh(\pi y)}$ \\
		\rule[-7.5pt]{0pt}{32pt} 
		8 & $\displaystyle \frac{\ds(2K' x,k)}{\cosh(\pi x)}$ & $\displaystyle  2\mi k\,\frac{\nd(2Ky,k')-e^{\pi y}}{e^{2\pi y}-1}$ 
		& 17 &
		$\displaystyle \frac{\dn(2K' x,k )}{\sinh(\pi x)}$ & $\displaystyle  -\mi\,\frac{\dn(2Ky,k')-e^{-\pi y}}{\sinh(\pi y)}$ \\
		\rule[-7.5pt]{0pt}{32pt} 
		9 & $\displaystyle \frac{\ns(2K' x,k)}{\cosh(\pi x)}$ & $\displaystyle  -2\mi k\,\frac{\nd(2Ky,k')}{e^{2\pi y}+1}$ 
		& 18 &
		$\displaystyle \frac{\sn(2K' x,k)}{\sinh(\pi x)}$ & $\displaystyle   \frac{\nd(2Ky,k')}{\cosh(\pi y)}$ 	 \\[10.0pt]
		\hline
	\end{tabular}
        \caption{Fourier transform pairs 1-18}\label{table:1}
\end{table}

To present our results on Fourier transform pairs, we define the Fourier transform operator in the following form:
\[
\mathcal{F}[f](y)=\int_{\RR+\i\epsilon} f(x)e^{2\pi\i xy}\d x,\quad y\in\RR,
\]
where $\epsilon$ is a small positive number. The inclusion of $\epsilon$ in this definition is necessary because the integrand may have poles on the real line. In all our examples, we work with functions $f$ that are analytic in a horizontal strip $0<\im(x)<b$ and decay exponentially fast as $\re(x) \to \pm \infty$. These conditions ensure that $\mathcal{F}[f](y)$ is well defined (the integral converges for $y\in \r$, and the result is independent of $\epsilon \in (0,b)$). 

In what follows, we assume that $k \in (0,1)$, and that $K=K(k)$ and $E=E(k)$ denote the complete elliptic integrals of the first and second kind, respectively, with $k'=\sqrt{1-k^2}$ and $K'=K(k')$ (see \cite[\S 19.2]{NIST}). We consider functions $f(x)$ of the form 
$J(2K'x,k)/\sinh(\pi x)$ or $J(2K'x,k)/\cosh(\pi x)$, where $J(\cdot, k)$ denotes one of the twelve Jacobi elliptic functions \cite[\S 22.2]{NIST}. With this scaling of variables, the function $f$ is periodic with period $2\i$ (see \cite[Table 22.4.3]{NIST}). Fourier transforms of these 24 functions are presented in Tables \ref{table:1} and \ref{table:2}. The results in Table \ref{table:1} correspond to functions $f$
that have only simple poles. In this case, ${\mathcal F}[f]$ remains within the same class of functions that can be represented as a ratio of a Jacobi elliptic function and $\sinh(\pi x)$ or $\cosh(\pi x)$ (up to multiplication by $\exp(\pm \pi y)$ and the addition of certain hyperbolic functions). The six identities in Table \ref{table:2} are more complicated, and the invariance observed above appears to be lost. These six identities correspond to cases in which $f$ has a double pole -- that is, when a zero of $\sinh(\pi x)$ or $\cosh(\pi x)$ coincides with a pole of the Jacobi elliptic function $J(2K'x,k)$.

\begin{table}[t]
    \centering
    \vspace{1mm}
	\begin{tabular}{|c|c|c|}
		\hline
		& $f(x)$ & ${\mathcal F}[f](y)$    \\
		\hline
		\rule[-7.5pt]{0pt}{35pt}
		19& $\displaystyle \frac{\cn(2K' x,k )}{\cosh(\pi x)}$ & $\displaystyle \frac{1}{k}\,\frac{1}{\sinh(\pi y)}\left\{k^2\int_0^{2Ky}\nd(u,k')^2\md u+2(E-K)y\right\}$\\
		\rule[-7.5pt]{0pt}{35pt}
		20& $\displaystyle \frac{\dn(2K' x,k )}{\cosh(\pi x)}$ & $\displaystyle \frac{1}{\sinh(\pi y)}\left\{\int_0^{2Ky}\dn(u,k')^2\md u+2(E-K)y\right\}$\\
		\rule[-7.5pt]{0pt}{35pt} 
		21& $\displaystyle \frac{\sn(2K' x,k )}{\cosh(\pi x)}$ & $\displaystyle \frac{\mi}{k}\,\frac{1}{\cosh(\pi y)}\left\{k^2\int_0^{2Ky}\nd(u,k')^2\md u+2(E-K)y\right\}$\\
		\rule[-7.5pt]{0pt}{35pt}
		22& $\displaystyle \frac{\cs(2K' x,k )}{\sinh(\pi x)}$ & $\displaystyle \frac{2}{1-e^{2\pi y}}\left\{\int_0^{2Ky}\dn(u,k')^2\md u+2(E-K)y\right\}$\\[15.0pt]		
		\rule[-7.5pt]{0pt}{35pt}
		23& $\displaystyle \frac{\ds(2K' x,k )}{\sinh(\pi x)}$ & $\displaystyle \frac{2}{1-e^{2\pi y}}\left\{k^2\int_0^{2Ky}\nd(u,k')^2\md u+2(E-K)y\right\}$\\
		\rule[-7.5pt]{0pt}{35pt}
		24& $\displaystyle \frac{\ns(2K' x,k)}{\sinh(\pi x)}$ & $\displaystyle \frac{2}{1+e^{2\pi y}}\left\{k^2\int_0^{2Ky}\nd(u,k')^2\md u+2(E-K)y\right\}$\\
		\hline
	\end{tabular}
        \caption{Fourier transform pairs 19-24}\label{table:2}
\end{table}

Our next results concern the Mellin transforms of functions $f(x)$ of the form $J(2K'x,k)/\sinh(\pi x)$ or $J(2K'x,k)/\cosh(\pi x)$ discussed above. By studying these Mellin transforms, we discovered sixteen functions of the form 
\begin{equation}\label{def_zeta_general2}
	\zeta_{j,l}(s,\tau):=\sideset{}{'}\sum_{\substack{n\in\ZZ\\m\geqslant 0}}(2a-1)^m (1-2b)^n\frac{\big [m+\frac{1}{2}+(n+\frac{d}{2})\tau \big]^{-s}}{\sin(\pi(\frac{c}{2}+(n+\frac{d}{2})\tau))},
\end{equation}
where $a,b,c,d\in\{0,1\}$, $j:=2a+b$, $l:=2c+d$, $\re(s)>1$, and $\im(\tau)>0$. The prime in the summation indicates that the term $n=0$ is omitted in the case $c=d=0$ (to avoid division by zero). Note that $j,l \in \{0,1,2,3\}$, and that the numbers $a,b$ (respectively, $c,d$) are the digits in the binary representation of $j$ (respectively, $l$). The reason for this seemingly complicated parameterization is that it reveals a certain symmetry between these functions under the transformation $(s,\tau) \mapsto (1-s,-1/\tau)$. 

In a series of papers \cite{Tsumura2008,TSUMURA_2009,Tsumura2010}, Tsumura obtained explicit expressions for many examples of Eisenstein-type series of the form
$$
\sum\limits_{\substack{m,n \in {\mathbb Z}\\ n\neq 0}}
(-1)^{n} \frac{[m+n\i]^{-j}}{h(\pi n)}, \;\;\; 
\sum\limits_{m,n \in {\mathbb Z}} 
\frac{[m+\frac{1}{2}+(n+\frac{1}{2})\i]^{-j}}{h(\pi (n+\frac{1}{2})) }, 
$$
where $j \in {\mathbb N}$ and the function $h$ is one of $\sinh$, $\cosh$, $\tanh$, or $\coth$. Some of these series are special cases (when $\tau=\i$) of our series defined in \eqref{def_zeta_general2}. Based on this analogy, we refer to the functions $\zeta_{j,l}(s,\tau)$ as the {\it Eisenstein-type series}. 

As in the case of the Riemann zeta and Dirichlet beta functions, we define the completed versions of the Eisenstein-type series as
\begin{equation}\label{def_completed_sum}
   \Lambda_{j,l}(s,\tau):= \pi^{-\frac{s}{2}}\Gamma\Big(\frac{s+\mathfrak{a}_{j,l}}{2}\Big)\zeta_{j,l}(s,\tau),
\end{equation}
where $\mathfrak{a}_{j,l}:=\lfloor (j+l)/2 \rfloor \; (\textnormal{mod } 2)$. 
As the second main result of our paper, we show that each of the sixteen functions $s\mapsto \Lambda_{j,l}(s,\tau)$ can be analytically continued to a meromorphic function in $\CC$, whose only singularity is a possible simple pole at $s=1$. When ${\mathfrak a}_{j,l}=0$ (respectively, ${\mathfrak a}_{j,l}=1$), the values of $\Lambda_{j,l}(n,\tau)$ are given explicitly (in terms of the elliptic integrals $K$, $K'$, and $E$) for positive even (respectively, odd) integers $n$. Moreover, each of these sixteen functions satisfies a functional equation that relates $\Lambda_{j,l}(s,\tau)$ and $\Lambda_{l,j}(1-s,-1/\tau)$. Interestingly, for fifteen of these functions $\Lambda_{j,l}(s,\tau)$, the results are essentially corollaries of the Fourier transform pairs presented in Tables \ref{table:1} and \ref{table:2}. However, the case of $\Lambda_{0,0}(s,\tau)$ is different, as it requires a Fourier transform pair (referred to as ``identity 25'') that is not included in Tables \ref{table:1} and \ref{table:2}. 

The paper is organized as follows. In Section \ref{section_Fourier_pairs}, we derive the 24 Fourier transform pairs given in Tables \ref{table:1} and \ref{table:2}, as well as Fourier transform identity 25, which is related to the logarithmic derivative of the Jacobi theta function $\theta_4(\cdot | \tau)$. The results concerning the Eisenstein-type series are presented and proved in Section \ref{section_Mellin}. In Section \ref{section_conclusion}, we offer several concluding remarks and discuss open problems related to the Eisenstein-type series.

\section{Fourier transform pairs}\label{section_Fourier_pairs}


First, we discuss how to prove the identities listed in Table \ref{table:1}. We note that identity 18 follows directly from identity 5 by applying the inverse Fourier transform. With a little additional work, one can verify that the same relationship holds for the pairs of identities $6\leftrightarrow 13$ and $14\leftrightarrow 16$. Using the periodicity properties of Jacobi elliptic functions \cite[Table 22.4.3]{NIST}, we observe that when $x$ is replaced by $x+\i/2$ or $x-\i/2$, the functions in Table \ref{table:1} are transformed (up to multiplicative constants) as follows:
 
	\begin{equation*}
		\begin{split}
			&\frac{\dc(2K'x,k)}{\cosh(\pi x)}\leftrightarrow \frac{\cd(2K'x,k)}{\sinh(\pi x)}\quad (1\leftrightarrow 13),\qquad \frac{\dc(2K'x,k)}{\sinh(\pi x)}\leftrightarrow \frac{\cd(2K'x,k)}{\cosh(\pi x)}\quad (10\leftrightarrow 4),\\[.5ex]
			&\frac{\nc(2K'x,k)}{\cosh(\pi x)}\leftrightarrow \frac{\sd(2K'x,k)}{\sinh(\pi x)}\quad (2\leftrightarrow 15),\qquad \frac{\nc(2K'x,k)}{\sinh(\pi x)}\leftrightarrow \frac{\sd(2K'x,k)}{\cosh(\pi x)}\quad (11\leftrightarrow 6),\\[.5ex]
			&\frac{\sc(2K'x,k)}{\cosh(\pi x)}\leftrightarrow \frac{\nd(2K'x,k)}{\sinh(\pi x)}\quad (3\leftrightarrow 14),\qquad \frac{\sc(2K'x,k)}{\sinh(\pi x)}\leftrightarrow \frac{\nd(2K'x,k)}{\cosh(\pi x)}\quad (12\leftrightarrow 5),\\[.5ex]
			&\frac{\cs(2K'x,k)}{\cosh(\pi x)}\leftrightarrow \frac{\dn(2K'x,k)}{\sinh(\pi x)}\quad (7\leftrightarrow 17),\qquad \frac{\ns(2K'x,k)}{\sinh(\pi x)}\leftrightarrow \frac{\sn(2K'x,k)}{\cosh(\pi x)}\quad (9\leftrightarrow 18),\\[.5ex]
			&\frac{\ds(2K'x,k)}{\cosh(\pi x)}\leftrightarrow \frac{\cn(2K'x,k)}{\sinh(\pi x)}\quad (8\leftrightarrow 16).
		\end{split}
	\end{equation*}
Thus, all identities in Table \ref{table:1} separate into the following six equivalence classes:
\begin{align*}
	1\leftrightarrow 13\leftrightarrow 6\leftrightarrow 11,&\quad 2\leftrightarrow 15,\quad 4\leftrightarrow 10,\quad 7\leftrightarrow 17,\\
	3\leftrightarrow 14\leftrightarrow 16\leftrightarrow 8,&\quad 12\leftrightarrow 5\leftrightarrow 18\leftrightarrow 9.
\end{align*}

Every identity within a given class follows from another identity in the same class by applying the inverse Fourier transform or performing a change of variables $x\mapsto x\pm \i/2$. Below we provide a complete proof of identity 1 and demonstrate how to derive identities 13 and 6 (thus illustrating the chain $1\rightarrow 13 \rightarrow 6$). Proofs of all other identities in Table \ref{table:1} can be obtained by exactly the same steps and are therefore omitted. 

\vspace{0.25cm}
\noindent
{\bf Proof of identity 1.}
We recall the following facts from \cite[\S 22]{NIST} concerning the Jacobi elliptic functions $\dc(z,k)$ and $\sd(z,k)$:
\begin{itemize}
	\item[(i)] $\dc(z,k)=\dc(z+2\mi K',k)=-\dc(z+2K,k)$;
    \item[(ii)] $\dc(\mi K',k)=k$;
	\item[(iii)] $\dc(z,k)$ has a simple pole at $z=K$ with residue $-1$, and all other poles are obtained by translations by $2mK+2n \i K'$ for $m,n\in {\mathbb Z}$;
	\item[(iv)] for $z\in\RR$ we have
	\begin{equation}\label{hyperbolic series for sd}
		\sd(z,k)=\frac{2\pi}{kk'K}\sum_{n\geqslant 0} (-1)^n\frac{q^{n+\frac{1}{2}}}{1+q^{2n+1}}\sin\big(\pi(n+\tfrac{1}{2})z/K\big),
	\end{equation}
    where $q:=\exp(-\pi K'/K)$. 
\end{itemize}

Fix $\epsilon\in (0,\frac{1}{2})$, $y\in\RR$, and $m\in {\mathbb N}$. Denote 
\begin{equation}\label{dc/cosh contour integral}
	I_m(y):=\int_{\gamma_m}\frac{\dc(2K'x,k)}{\cosh(\pi x)}e^{2\pi\mi xy}\md x,
\end{equation}
where $R>0$, and the contour $\gamma_m$ is the boundary of the rectangle with vertices 
\[
\bigl\{-mK/K'+\mi\epsilon,\, mK/K'+\mi\epsilon,\, mK/K'+\mi(\epsilon+1),\, -mK/K'+\mi(\epsilon+1)\bigr\},\]
traversed in the counterclockwise direction. From the facts listed in items (i)-(iii) above, it follows that the function $x\mapsto \dc(2K'x,k)$ is analytic in the strip $0<\im(x)<3/2$, except for simple poles at the points $x_n:=\mi+(n+1/2)K/K'$, $n\in\ZZ$, with corresponding residues $\left(-1\right)^{n+1}/(2K')$. The function $1/\cosh(\pi x)$ has a simple pole at $x=\i/2$ with residue $-\i/\pi$. 

We apply the Cauchy Residue Theorem to evaluate the integral in \eqref{dc/cosh contour integral} as follows:
\begin{align*}
	I_m(y)& =2e^{-\pi y}\,\dc(\mi K',k)+2\pi\mi\sum_{n=-m}^{m-1}\frac{(-1)^{n+1}}{2K'}\frac{\exp(2\pi\mi x_n y)}{\cosh(\pi x_n)}\\
	& =2k\,e^{-\pi y}+\frac{\pi\mi}{K'}e^{-2\pi y}\sum_{n=-m}^{m-1}\left(-1\right)^n\frac{\exp\left((2n+1)\pi\mi\frac{K}{K'}y\right)}{\cosh\left((n+\frac{1}{2})\pi\frac{K}{K'}\right)}\\
    & =2k\,e^{-\pi y}-\frac{2\pi}{K'}e^{-2\pi y}\sum_{n=0}^{m-1}\left(-1\right)^n\frac{\sin\left((2n+1)\pi\frac{K}{K'}y\right)}{\cosh\left((n+\frac{1}{2})\pi\frac{K}{K'}\right)}.
\end{align*}

Next, we let $m \to +\infty$. The integrals over the vertical sides of $\gamma_m$ vanish as $m\to+\infty$. To justify this, note that $|\dc(2K'x\pm 2mK,k)|$ does not depend on $m$ (see item (i) above), and the function $\dc(2K'x,k)$ has no poles on the vertical sides of $\gamma_m$, while $1/\cosh(\pi x)\to 0$ as $|\re(x)|\to \infty$. \label{limit_m_infinity} Thus, in the limit as $m\to+\infty$, we obtain
\begin{align*}
	\int_{\RR+\mi\epsilon}\frac{\dc(2K'x,k)}{\cosh(\pi x)}&e^{2\pi\mi xy}\md x-\int_{\RR+\mi(\epsilon+1)}\frac{\dc(2K'x,k)}{\cosh(\pi x)}e^{2\pi\mi xy}\md x\\
	&=2k\,e^{-\pi y}-\frac{2\pi}{K'}e^{-2\pi y}\sum_{n\geqslant 0}\left(-1\right)^n\frac{\sin\left((2n+1)\pi\frac{K}{K'}y\right)}{\cosh\left((n+\frac{1}{2})\pi\frac{K}{K'}\right)}.
\end{align*}
Changing the variable of integration $x\mapsto x+\mi$ in the second integral and using the periodicity property of $\dc(z,k)$ stated in item (i) above, we obtain
\begin{align*}
	\left(1+e^{-2\pi y}\right)&\int_{\RR+\mi\epsilon}\frac{\dc(2K'x,k)}{\cosh(\pi x)}e^{2\pi\mi xy}\md x=2k\,e^{-\pi y}-\frac{2\pi}{K'}e^{-2\pi y}\sum_{n\ge 0}\left(-1\right)^n\frac{\sin\left((2n+1)\pi\frac{K}{K'}y\right)}{\cosh\left((n+\frac{1}{2})\pi\frac{K}{K'}\right)}.
\end{align*}
The infinite series can be evaluated with the help of formula \eqref{hyperbolic series for sd}, by first replacing $k \mapsto k'$ and then $z\mapsto 2 K y$:
\[
\frac{2\pi}{K'}\sum_{n=0}^{\infty}\left(-1\right)^n\frac{\sin\left((2n+1)\pi\frac{K}{K'}y\right)}{\cosh\left((n+\frac{1}{2})\pi\frac{K}{K'}\right)}=2kk'\sd(2Ky,k').
\]
The above two equations together imply identity~1.
\qed 

\vspace{0.25cm}
\noindent
{\bf Deriving identity 13 from identity 1.}
\label{deriving_13_from_1}
Take $\epsilon\in (0,\frac{1}{2})$ and denote $\delta:=\frac{1}{2}-\epsilon$ and
\[
I(y):=\int_{\RR+\mi\epsilon}\frac{\cd(2K'x,k)}{\sinh(\pi x)}e^{2\pi\mi xy}\md x,\quad y\in\RR.
\]
Note that the function $\cd(z,k)$ has simple poles at $(2m+1)K+(2n+1)\mi K'$, $m,n\in\ZZ$, and satisfies $\cd(0,k)=1$ (see \cite[Table 22.4.1]{NIST} and \cite[Table 22.5.1]{NIST}).  
Replace the path of integration $\RR+\mi\epsilon$ by the positively oriented contour $\gamma_m$, which is the boundary of a rectangle with vertices 
\[
\bigl\{-mK/K'-\mi\delta,\, mK/K'-\mi\delta,\, mK/K'+\mi\epsilon,\, -mK/K'+\mi\epsilon\bigr\},\]
 and let $m\to+\infty$.  
Applying the Cauchy Residue Theorem (and eliminating the integrals over the vertical sides of $\gamma_m$ by an argument similar to that used in the proof of identity 1 on page~\pageref{limit_m_infinity}) gives
\[
\int_{\RR-\mi \delta}\frac{\cd(2K'x,k)}{\sinh(\pi x)}e^{2\pi\mi xy}\md x-\int_{\RR+\mi\epsilon}\frac{\cd(2K'x,k)}{\sinh(\pi x)}e^{2\pi\mi xy}\md x
=2\pi\mi \times \Res_{x=0}\frac{\cd(2K'x,k)}{\sinh(\pi x)}e^{2\pi\mi xy}.
\]
Setting $x=u-\i/2$ in the first integral gives
\[
I(y)=-2\mi+e^{\pi y} \int_{\RR+\mi\epsilon}\frac{\cd(2K'u-\mi K',k)}{\sinh(\pi u-\pi\mi/2)}e^{2\pi\mi uy}\md u.
\]
According to \cite[Table 22.4.3]{NIST}, we have
\[
\cd(z+2\mi K',k)=\cd(z,k),\quad \cd(z+\mi K',k)=k^{-1}\,\dc(z,k),
\]
therefore,
\[
I(y)=-2\mi+\frac{\mi}{k}e^{\pi y}\int_{\RR+\mi\epsilon}\frac{\dc(2K'x,k)}{\cosh(\pi x)}e^{2\pi\mi xy}\md x
=-\mi\,\frac{k'\,\sd(2Ky,k')+e^{-\pi y}}{\cosh(\pi y)},
\]
where, in the last step, we applied identity~1 from Table~\ref{table:1}.
\qed

\vspace{0.25cm}
\noindent
{\bf Deriving identity 6 from identity 13.}
Fix $\epsilon\in(0,\frac{1}{2})$ and let
\[
f(x):=\frac{\cd(2Kx,k')-1}{\sinh(\pi x)}.
\]
Then $f(x)$ is analytic in the strip $-\frac{1}{2}<\im(x)<\frac{1}{2}$ and is a Schwartz function. We will also need the following result:
\begin{equation}\label{Fourier/sinh}
\int_{\RR+\mi\epsilon}\frac{e^{2\pi\mi xy}}{\sinh(\pi x)}\md x=-\mi\,\frac{e^{-\pi y}}{\cosh(\pi y)}.
\end{equation}
The simplest way to establish \eqref{Fourier/sinh} is to use the identity $2\cos(xy)=\exp(\i xy)+\exp(-\i xy)$ to express the integral in \eqref{Fc_cosh} as an integral over $\r$, then change the variable of integration $x\mapsto x-\pi \i/2$ and rescale the variables $x$ and $y$.

Using \eqref{Fourier/sinh} and identity 13 in Table~\ref{table:1}, we obtain for $y\in \r$
\begin{align*}
	g(y)& :=\int_{\RR}f(x)e^{2\pi\mi xy}\md x=\int_{\RR+\mi\epsilon}f(x)e^{2\pi\mi xy}\md x\\
	& =-\mi\,\frac{k\, \sd\left(2K'y;k\right)+e^{-\pi y}}{\cosh(\pi y)}+\mi\,\frac{e^{-\pi y}}{\cosh(\pi y)} =-\mi k\,\frac{\sd\left(2K'y;k\right)}{\cosh(\pi y)}.
\end{align*}
Writing $f$ as the inverse Fourier transform of $g$ (see \cite[Theorem 2.2.14]{Graf}), we have 
\[
-\mi k\int_{\RR}\frac{\sd\left(2K'x;k\right)}{\cosh(\pi x)}e^{-2\pi\mi xy}\md x
=\frac{\cd(2Ky,k')-1}{\sinh(\pi y)},\;\;\; y\in \r,
\]
which is equivalent to identity~6 in Table~\ref{table:1}.
\qed

Next, we turn our attention to proving the identities listed in Table~\ref{table:2}. Using \cite[\S 22.4]{NIST}, one can verify that the functions $f$ appearing in Table~\ref{table:1} have only simple poles, whereas each function $f$ appearing in identities~19--21 (respectively, 22--24) in Table~\ref{table:2} has a double pole at $x=\i/2$ (respectively, $x=0$).  
Moreover, under the shift $x\mapsto x \pm \i/2$, the following transformation relations hold (up to a multiplicative constant):
\begin{equation*}
	\begin{split}
		\frac{\sn(2K'x,k)}{\cosh(\pi x)}\leftrightarrow \frac{\ns(2K'x,k)}{\sinh(\pi x)},&\quad (21\leftrightarrow 24),\\[.5ex]
		\frac{\cn(2K'x,k)}{\cosh(\pi x)}\leftrightarrow \frac{\ds(2K'x,k)}{\sinh(\pi x)},&\quad (19\leftrightarrow 23),\\[.5ex]
		\frac{\dn(2K'x,k)}{\cosh(\pi x)}\leftrightarrow \frac{\cs(2K'x,k)}{\sinh(\pi x)},&\quad (20\leftrightarrow 22).
	\end{split}
\end{equation*}
We will prove identity~21. Using the same technique as in the derivation of identity~13 from identity~1 on page~\pageref{deriving_13_from_1}, one can then deduce identity~24. The proofs of identities~19, 20, 22, and~23 follow exactly the same steps and are therefore omitted.

\vspace{0.25cm}
\noindent
{\bf Proof of identity 21.}
For $y\in\RR$ and $m\in {\mathbb N}$, define
\begin{equation}\label{def_I_m3}
I_m(y):=\int_{\gamma_m}\frac{\sn(2K'x,k)}{\cosh(\pi x)}e^{2\pi\mi xy}\md x,
\end{equation}
where $\gamma_m$ is the boundary of the rectangle with vertices 
\[
\bigl\{-(m+\tfrac{1}{2})K/K'+\mi\epsilon,\,
  (m+\tfrac{1}{2})K/K'+\mi\epsilon,\,
  (m+\tfrac{1}{2})K/K'+\mi(\epsilon+1),\,
  -(m+\tfrac{1}{2})K/K'+\mi(\epsilon+1)\bigr\},
\]
traversed in the counterclockwise direction.  
Recall the Laurent expansion of $\sn(x,k)$ at $x=\i K'$ (see \cite[Eq.~(2.43)]{ArEb}):
\[
\sn(u+\i K',k)=\frac{1}{ku}+\frac{1+k^2}{6k}u+\mo\!\left(u^3\right), \quad u\to 0.
\]
The integrand in \eqref{def_I_m3} has a double pole at $x=\i/2$ with residue
\[
\Res_{x=\mi/2}\frac{\sn(2K'x,k)}{\cosh(\pi x)}e^{2\pi\mi xy}
=\frac{y\,e^{-\pi y}}{kK'}.
\]
Moreover, the integrand has simple poles within $\gamma_m$ at the points 
$x_n:=\i/2+nK/K'$, $0<|n|\le m$, with residues $\left(-1\right)^n/(2kK')$.  
Applying the Cauchy Residue Theorem yields
\begin{align*}
	I_m(y)& =2\pi\mi\frac{y\,e^{-\pi y}}{kK'}
    +2\pi\mi \sideset{}{'}\sum_{n=-m}^m\frac{(-1)^n}{2kK'}\frac{\exp(2\pi\mi x_n y)}{\cosh(\pi x_n)}\\
	& =\frac{2\pi\mi y e^{-\pi y}}{kK'}
    +\frac{\pi e^{-\pi y}}{kK'}\sideset{}{'}\sum_{n=-m}^m(-1)^n
      \frac{\exp(2n\pi\mi\frac{K}{K'}y)}{\sinh(n\pi\frac{K}{K'})},
\end{align*}
where the prime in the summation indicates that the term with $n=0$ is omitted.  
Letting $m\to+\infty$ and combining the terms with indices $\pm n$, we obtain
\begin{align*}
\int_{\RR}\frac{\sn(2K'x,k)}{\cosh(\pi x)}&e^{2\pi\mi xy}\md x
-\int_{\RR+\mi}\frac{\sn(2K'x,k)}{\cosh(\pi x)}e^{2\pi\mi xy}\md x\\
&=\frac{2\pi\mi y\,e^{-\pi y}}{kK'}
  +\frac{2\pi\mi\,e^{-\pi y}}{kK'}
     \sum_{n \geqslant 1}(-1)^n
     \frac{\sin(2n\pi\frac{K}{K'}y)}{\sinh(n\pi\frac{K}{K'})}.
\end{align*}
Changing the variable of integration $x\mapsto x+\mi$ in the second integral, we arrive at
\begin{equation}\label{sn/cosh Fourier transform-1}
\int_{\RR}\frac{\sn(2K'x,k)}{\cosh(\pi x)}e^{2\pi\mi xy}\md x
=\frac{\pi\mi}{kK'\cosh(\pi y)}
 \bigg\{y+\sum_{n \geqslant 1}(-1)^n
 \frac{\sin(2n\pi\frac{K}{K'}y)}{\sinh(n\pi\frac{K}{K'})}\bigg\},
\end{equation}
where we used the fact that $\sn(z+2\mi K',k)=\sn(z,k)$ 
(see \cite[Table~22.4.3]{NIST}).

Recall the Fourier series \cite[p.~25, Eq.~(2.23)]{Ober}:
\[
\sum_{n \geqslant 1}\frac{nq^n}{1-q^{2n}}\cos(2n\pi z)
=\frac{K^2}{2\pi^2}\dn(2Kz,k)^2-\frac{KE}{2\pi^2}, \;\;\; z\in \r,
\]
where $q=\exp(-\pi K'/K)$.  
Replacing $z$ by $z+1/2$ and using the identity 
$\dn(z+K,k)=k' \nd(z,k)$ (see \cite[Table~22.4.3]{NIST}), we find
\[
\sum_{n \geqslant 1}(-1)^n\frac{nq^n}{1-q^{2n}}\cos(2n\pi z)
=\frac{K^2 k'^2}{2\pi^2}\nd(2Kz,k)^2-\frac{KE}{2\pi^2}.
\]
Next, we replace $k$ by $k'$ and set $z=uK/K'$, which gives
\[
\sum_{n \geqslant 1}\left(-1\right)^n
\frac{n\,\cos\!\left(2n\pi\frac{K}{K'}u\right)}{\sinh\!\left(n\pi\frac{K}{K'}\right)}
=\frac{K'^2 k^2}{\pi^2}\nd\!\left(2Ku;k'\right)^2-\frac{K'E'}{\pi^2},
\]
where $E':=E(k')$.  
Integrating both sides over the interval $[0,y]$ gives
\[
\sum_{n \geqslant 1}(-1)^n 
 \frac{\sin(2n\pi\frac{K}{K'}y)}{\sinh(n\pi\frac{K}{K'})}
 =\frac{2K}{\pi}\left\{K'k^2\int_0^y\nd(2Ku,k')^2\md u-E'y\right\}.
\]
Combining this with \eqref{sn/cosh Fourier transform-1} yields
\[
\int_{\RR}\frac{\sn(2K'x,k)}{\cosh(\pi x)}e^{2\pi\mi xy}\md x
=\frac{\mi}{\cosh(\pi y)}
 \left\{k\int_0^{2Ky}\nd(u,k')^2\md u
 +\frac{\pi-2KE'}{kK'}y\right\},
\]
which is equivalent to identity~21, in view of Legendre’s relation for complete elliptic integrals \cite[Eq.~(19.7.1)]{NIST}.
\qed

We also record here the following result about the Fourier transform pair related to the function 
  \begin{equation}\label{def_mathfrak_f}
	\mathfrak{f}(x,k):=\frac{\theta_4'(-\mi\pi\tau x|\tau)}{\theta_4(-\mi\pi\tau x|\tau)},
	\end{equation}
	where $\tau:=\i K'/K$. This result will be needed in the next section, where we will use it to study the Eisenstein-type series $\zeta_{0,0}(s,\tau)$ defined by 
	\eqref{def_zeta_general2}.

\vspace{0.25cm}
\noindent
{\bf Identity 25:}
\label{prop_f_cosh_Fourier}
\begin{equation}\label{g(x,k) Fourier transform}
	\int_{\RR}\frac{\mathfrak{f}(x,k)}{\cosh(\pi x)}e^{2\pi\mi xy}\md x=\frac{\mi K}{K'}\frac{1}{ \cosh(\pi y)}\left\{\mathfrak{f}(y,k')+2y-\frac{K'}{K}\tanh(\pi y)\right\}, \;\;\; y\in \r.
\end{equation}
\begin{proof}
The proof is very similar to the proof of identity 1 that we presented above; therefore, we only provide a sketch of the main steps of the proof.
Using the properties of theta functions in \cite[\S 20]{NIST}, we check that ${\mathfrak f}(x,k)$ satisfies the periodicity properties 
$$
\mathfrak{f}(x+\mi,k)=\mathfrak{f}(x,k)-2\i, 
\;\;\;
\mathfrak{f}(x+K/K',k)=\mathfrak{f}(x,k). 
$$
The function ${\mathfrak f}(\cdot,k)$ has simple poles at $x_n:=\i/2+nK/K'$, $(n\in\ZZ)$ with residues $\i /(\pi \tau)$, 
 and a Fourier series expansion
\begin{equation}\label{prop1_proof1}
	\mathfrak{f}(x,k)=2 \i \sum\limits_{n\geqslant 1}
	\frac{\sin(2\pi n \frac{K'}{K} x)}{\sin(\pi n \tau)}, \;\;\; 
	x\in \r. 
\end{equation} 

Denote $g(x,k):={\mathfrak f}(x,k)/\cosh(\pi x)$. The function $g$ has a double pole at $x=x_0=\i/2$ and simple poles at points $x_n$, $n \in {\mathbb Z} \setminus \{0\}$.   
Shifting the contour of integration we have 
\begin{equation}\label{prop1_proof2}
\int_{\RR}g(x,k)e^{2\pi\mi xy}\md x-\int_{\RR+\mi}g(x,k)e^{2\pi\mi xy}\md x=2\pi\mi\sum_{n\in\ZZ}\Res_{x=x_n} g(x,k)e^{2\pi\mi xy}.
\end{equation}
Applying \cite[Eq. (20.2.11)]{NIST} and \cite[Eq. (20.6.2)]{NIST}, we find that as $z\to 0$,
\[\theta_4(z+\pi \tau/2 |\tau)=\mi e^{-\mi z-\frac{\pi}{4}\mi\tau}\theta_1(z|\tau)=\mi e^{-\frac{\pi}{4}\mi\tau}\theta_1'(0|\tau)z\big(1-\mi z+O(z^2)\big).\]
With the help of the above result, we compute the residue at the double pole $x=\i/2$ 
\[\Res_{x=\frac{\mi}{2}}g(x,k)e^{2\pi\mi xy}=\frac{e^{-\pi y}}{\pi\tau}(2\mi y-\tau),\]
and the residues at the simple poles are given by
\[\Res_{x=x_n}g(x,k)e^{2\pi\mi xy}=
-\frac{\mi e^{-\pi y}}{\pi\tau}\frac{e^{-2\pi n y/\tau}}{\sin(\pi n/\tau)}, \;\;\; n \in {\mathbb Z} \setminus \{0\}. \]
Combining the above two results with \eqref{prop1_proof1} and 
\eqref{prop1_proof2} leads to 
\eqref{g(x,k) Fourier transform}.
\end{proof}

\section{Eisenstein-type series}
\label{section_Mellin}

In this section, we state and prove results concerning analyticity, functional equations, and values at positive odd or even integers for the sixteen Eisenstein-type series defined in \eqref{def_zeta_general2}. We first state these results in a sequence of ten propositions, and then discuss their proofs.  

In the previous section, we focused on Jacobi elliptic functions with a real parameter $k\in (0,1)$. In this section, it is more convenient to parameterize everything in terms of the parameter $\tau$ lying in the upper half-plane $\im(\tau)>0$. Then, from \cite[\S 20.9]{NIST}, we find 
\begin{equation}\label{from_tau_to_k}
k=\frac{\theta_2^2(0|\tau)}{\theta_3^2(0|\tau)},
\;\;\; k'=\frac{\theta_4^2(0|\tau)}{\theta_3^2(0|\tau)},
\;\;\; K=\frac{\pi}{2} \theta_3^2(0| \tau), \;\;\;
 K'=-\i \tau K. 
\end{equation} 
Note that changing $k \mapsto k'$ is equivalent to transforming $\tau \mapsto \tau':=-1/\tau$. We also remind the reader that the functions $\Lambda(s)$ and $\Lambda(s,\chi_4)$ are defined in equations \eqref{def_Lambda_zeta} and \eqref{def_Lambda_beta}.

\begin{proposition}\label{Prop_1}
	For $\im(\tau)>0$ and $\re(s)>0$ define 
	\begin{equation*}
		\Lambda_{1,1}(s,\tau):=\pi^{-\frac{s}{2}}\Gamma\Big(\frac{1+s}{2}\Big)\sum_{\substack{n\in\ZZ\\m\geqslant 0}} (-1)^{m+n}\frac{\big[m+\frac{1}{2}+(n+\frac{1}{2})\tau\big]^{-s}}{\sin(\pi (n+\frac{1}{2} )\tau)}.
	\end{equation*}
The function $\Lambda_{1,1}(s,\tau)$ can be analytically continued to an entire function of $s$ and it satisfies 
	\begin{equation*}
		\Lambda_{1,1}(s,\tau)+\frac{\mi}{\tau}\Lambda_{1,1}\big(1-s,-1/\tau\big)=-\mi\,\theta_3^2(0|\tau)\Lambda(s,\chi_4), \;\;\; s\in \c.  
	\end{equation*}
	Moreover, for all $n\in\ZZ_{\geqslant 0}$,
	\begin{equation*}
		\Lambda_{1,1}(2n+1,\tau)=\mi K\frac{(-4\pi)^{-n}}{\Gamma(n+\frac{1}{2})}\left[\frac{\md^{2n}}{\md y^{2n}}\frac{k'\;\cd(2Ky,k')-1}{\cosh(\pi y)}\right]_{y=0}.
	\end{equation*}
	In particular,
	\[\Lambda_{1,1}(1,\tau)=\frac{\mi (k'-1)K}{\pi^{1/2}},\quad\Lambda_{1,1}(3,\tau)=\frac{\mi K}{2\pi^{3/2}}\left[\pi^2(k'-1)+4k'k^2 K^2\right].\]
\end{proposition}

\begin{proposition}\label{Prop_2}
For $\im(\tau)>0$ and $\re(s)>0$ define 
	\begin{align*}
		\Lambda_{1,2}(s,\tau)&:=\pi^{-\frac{s}{2}}\Gamma\Big(\frac{1+s}{2}\Big)\sum_{\substack{n\in\ZZ\\m\geqslant 0}}(-1)^{m+n}\frac{\big[m+\frac{1}{2}+n\tau\big]^{-s}}{\cos(\pi n\tau)},\\
		\Lambda_{2,1}(s,\tau)&:=\pi^{-\frac{s}{2}}\Gamma\Big(\frac{1+s}{2}\Big)\sum_{\substack{n\in\ZZ\\m\geqslant 0}}\frac{\big[m+\frac{1}{2}+(n+\frac{1}{2})\tau\big]^{-s}}{\sin(\pi(n+\frac{1}{2})\tau)}.
	\end{align*}
The functions $\Lambda_{1,2}(s,\tau)$ and $\Lambda_{2,1}(s,\tau)$  can be analytically continued to entire functions of $s$ and they satisfy
	\begin{equation}
		\Lambda_{1,2}(s,\tau)=-\frac{\mi}{\tau}\Lambda_{2,1}\big(1-s,-1/\tau\big), \;\;\; s\in \c. 
	\end{equation}
	Moreover, for all $n\in\ZZ_{\geqslant 0}$,
	\begin{align*}
		\Lambda_{1,2}(2n+1,\tau)& =k'K\frac{(-4\pi)^{-n}}{\Gamma(n+\frac{1}{2})}\left[\frac{\md^{2n}}{\md y^{2n}}\frac{\sn(2Ky,k')}{\sinh(\pi y)}\right]_{y=0}, \\[.5ex]
		\Lambda_{2,1}(2n+1,\tau)& =-kK \frac{(-4\pi)^{-n}}{\Gamma(n+\frac{1}{2})}\left[\frac{\md^{2n}}{\md y^{2n}}\frac{\nd(2Ky,k')}{\cosh(\pi y)}\right]_{y=0}. 
	\end{align*}
	In particular,
	\begin{alignat*}{2}
		\Lambda_{1,2}(1,\tau)& =\frac{2k'K^2}{\pi^{3/2}},&\quad\Lambda_{1,2}(3,\tau)&=\frac{k'K^2}{3\pi^{5/2}}\left[\pi^2+4(2-k^2)K^2\right],\\[.5ex]
		\Lambda_{2,1}(1,\tau)& =-\frac{kK}{\pi^{1/2}},&\quad\Lambda_{2,1}(3,\tau)&=-\frac{kK}{2\pi^{3/2}}\left[\pi^2+4(k^2-1)K^2\right].
	\end{alignat*}
\end{proposition}

\begin{proposition}\label{Prop_3}
For $\im(\tau)>0$ and $\re(s)>0$ define
	\begin{equation}\label{def_Lambda13}
		\Lambda_{1,3}(s,\tau):=\pi^{-\frac{s}{2}}\Gamma\Big(\frac{s}{2}\Big)\sum_{\substack{n\in\ZZ\\m\geqslant 0}}\left(-1\right)^{m+n}\frac{\big[m+\frac{1}{2}+(n+\frac{1}{2})\tau\big]^{-s}}{\cos(\pi (n+\frac{1}{2})\tau)},
	\end{equation}
	and for $\re(s)>1$ define
	\begin{equation}\label{def_Lambda31}
		\Lambda_{3,1}(s,\tau):=\pi^{-\frac{s}{2}}\Gamma\Big(\frac{s}{2}\Big)\sum_{\substack{n\in\ZZ\\m\geqslant 0}}(-1)^n\frac{\big[m+\frac{1}{2}+(n+\frac{1}{2})\tau\big]^{-s}}{\sin(\pi (n+\frac{1}{2} )\tau )}.
	\end{equation}
	The function $\Lambda_{1,3}(s,\tau)$ can be analytically continued to an enture function of $s$ and the function $\Lambda_{3,1}(s,\tau)$ can be extended to an analytic function in $\c \setminus \{1\}$, having a simple pole at $s=1$ with residue $-2\i k K/\pi$. These functions satisfy
	\begin{equation}\label{reflection formula for Lambda_1,3 and Lambda_3,1}
		\Lambda_{1,3}(s,\tau)+\frac{\mi}{\tau}\,\Lambda_{3,1}\big(1-s,-1/\tau\big)=-\mi\,\theta_4^2(0|\tau)(2^{1-s}-1)\Lambda(s), \;\;\; s\in \c. 
	\end{equation}
	Moreover, for all $n\in\ZZ_{\geqslant 0}$,
	\begin{align}
		\Lambda_{1,3}(2n+2,\tau)& =\frac{\mi k'K}{2\sqrt{\pi}}\frac{(-4\pi)^{-n}}{\Gamma(n+\frac{3}{2})}\left[\frac{\md^{2n+1}}{\md y^{2n+1}}\frac{\cd(2Ky,k')-1}{\sinh(\pi y)}\right]_{y=0}, \label{Lambda_1,3 values}\\[.5ex]
		\Lambda_{3,1}(2n+2,\tau)& =\frac{\mi kK}{2\sqrt{\pi}}\frac{(-4 \pi)^{-n}}{\Gamma(n+\frac{3}{2})}\left[\frac{\md^{2n+1}}{\md y^{2n+1}}\frac{k'\;\sd(2Ky,k')+e^{-\pi y}}{\cosh(\pi y)}\right]_{y=0}. \label{Lambda_3,1 values}
	\end{align}
	In particular,
	\begin{alignat*}{2}
		\Lambda_{1,3}(2,\tau)&=\frac{2k'k^2 K^3}{\pi^2\mi},&\quad\Lambda_{1,3}(4,\tau)&=\frac{k'k^2 K^3}{3\pi^3\mi}\left[\pi^2+2(5k^2-4)K^2\right],\\[.5ex]
		\Lambda_{3,1}(2,\tau)&=\frac{kK}{\pi\mi}(\pi-2k'K), &\quad\Lambda_{3,1}(4,\tau)&=\frac{kK}{3\pi^2\mi}\left[\pi^3-3\pi^2 k'K+4k'(1-2k^2)K^3\right].
	\end{alignat*}
\end{proposition}

\begin{proposition}\label{Prop_4}
For $\im(\tau)>0$ and $\re(s)>1$ define 
	\begin{equation*}
		\Lambda_{2,2}(s,\tau):=\pi^{-\frac{s}{2}}\Gamma\Big(\frac{s}{2}\Big)\sum_{\substack{n\in\ZZ\\m\geqslant 0}}\frac{\big[ m+\frac{1}{2} +n\tau\big]^{-s}}{\cos(\pi n\tau )}.
	\end{equation*}
	The function $\Lambda_{2,2}(s,\tau)$ can be extended to an analytic function in $\c \setminus \{1\}$, having a simple pole at $s=1$ with residue $2K/\pi$, and satisfying 
	\begin{equation*}
		\Lambda_{2,2}(s,\tau)+\frac{\mi}{\tau}\,\Lambda_{2,2}\big(1-s,-1/\tau\big)=\theta_3^2\left(0|\tau\right)\Lambda(s), \;\;\; s\in \c. 
	\end{equation*}
	Moreover, for all $n\in\ZZ_{\geqslant 0}$,
	\begin{equation*}
		\Lambda_{2,2}(2n+2,\tau)=-\frac{K}{2\sqrt{\pi}}\frac{(-4 \pi)^{-n}}{\Gamma(n+\frac{3}{2})}\left[\frac{\md^{2n+1}}{\md y^{2n+1}}\frac{\dn(2Ky,k')-e^{-\pi y}}{\sinh(\pi y)}\right]_{y=0}.
	\end{equation*}
	In particular,
	\begin{align*}
		\Lambda_{2,2}(2,\tau)&=\frac{K}{2\pi^2}\left[\pi^2+4(1-k^2)K^2\right],\\[.5ex]
		\Lambda_{2,2}(4,\tau)&=\frac{K}{24\pi^3}\left[\pi^4+8\pi^2(1-k^2)K^2+16(5-6k^2+k^4)K^4\right].
	\end{align*}
\end{proposition}

\begin{proposition}\label{Prop_5}
	For $\im(\tau)>0$ and $\re(s)>1$ define 
	\begin{align*}
		\Lambda_{3,2}(s,\tau)&:=\pi^{-\frac{s}{2}}\Gamma\Big(\frac{s}{2}\Big)\sum_{\substack{n\in\ZZ\\m\geqslant 0}}(-1)^n\frac{\big[m+\frac{1}{2}+n\tau\big]^{-s}}{\cos(\pi n\tau)},\\
		\Lambda_{2,3}(s,\tau)&:=\pi^{-\frac{s}{2}}\Gamma\Big(\frac{s}{2}\Big)\sum_{\substack{n\in\ZZ\\m\geqslant 0}}\frac{\big[m+\frac{1}{2}+(n+\frac{1}{2})\tau\big]^{-s}}{\cos(\pi(n+\frac{1}{2})\tau)}.
	\end{align*}
	The functions  $\Lambda_{3,2}(s,\tau)$ and $\Lambda_{2,3}(s,\tau)$ can be extended to analytic functions in $\c \setminus \{1\}$, having simple pole at $s=1$ with the corresponding residues $2k'K/\pi$ and $2kK/\pi$. These functions  satisfy
	\begin{equation*}
		\Lambda_{3,2}(s,\tau)+\frac{\mi}{\tau}\,\Lambda_{2,3}\big(1-s,-1/\tau\big)=\theta_4^2\left(0|\tau\right)\Lambda(s), \;\;\; s\in \c. 
	\end{equation*}
	Moreover, for all $n\in\ZZ_{\geqslant 0}$,
	\begin{align*}
		\Lambda_{3,2}(2n+2,\tau)& =-\frac{k'K}{2\sqrt{\pi}}\frac{(-4 \pi)^{-n}}{\Gamma(n+\frac{3}{2})}\left[\frac{\md^{2n+1}}{\md y^{2n+1}}\frac{\cn(2Ky,k')-e^{-\pi y}}{\sinh(\pi y)}\right]_{y=0},\\[.5ex]
		\Lambda_{2,3}(2n+2,\tau)& =-\frac{kK}{2\sqrt{\pi}}\frac{(-4 \pi)^{-n}}{\Gamma(n+\frac{3}{2})}\left[\frac{\md^{2n+1}}{\md y^{2n+1}}\frac{\nd(2Ky,k')-e^{-\pi y}}{\sinh(\pi y)}\right]_{y=0}. 
	\end{align*}
	In particular,
	\[\Lambda_{3,2}(2,\tau)=\frac{k'K}{2\pi^2}\left(\pi^2+4K^2\right),\quad \Lambda_{2,3}(2,\tau)=\frac{kK}{2\pi^2}\left[\pi^2+4(k^2-1)K^2\right].\]
\end{proposition}

\begin{proposition}\label{Prop_6}
	For $\im(\tau)>0$ and $\re(s)>1$ define
	\begin{equation*}
		\Lambda_{3,3}(s,\tau):=\pi^{-\frac{s}{2}}\Gamma\Big(\frac{1+s}{2}\Big)\sum_{\substack{n\in\ZZ\\m\geqslant 0}}(-1)^n\frac{\big[m+\frac{1}{2}+(n+\frac{1}{2})\tau\big]^{-s}}{\cos(\pi(n+\frac{1}{2})\tau)}.
	\end{equation*}
	The function $\Lambda_{3,3}(s,\tau)$ can be analytically continued to an entire function of $s$ and it satisfies 
	\begin{equation*}
		\Lambda_{3,3}(s,\tau)=\frac{\mi}{\tau}\,\Lambda_{3,3}\big(1-s,-1/\tau\big), \;\;\; s\in \c. 
	\end{equation*}
	Moreover, for all $n\in\ZZ_{\geqslant 0}$,
	\begin{equation*}
		\Lambda_{3,3}(2n+1,\tau)=-\i k k' K \frac{(-4\pi)^{-n}}{\Gamma(n+\frac{1}{2})}\left[\frac{\md^{2n}}{\md y^{2n}}\frac{\sd(2Ky,k')}{\sinh(\pi y)}\right]_{y=0}.
	\end{equation*}
	In particular,
	\[\Lambda_{3,3}(1,\tau)=-\i \frac{2kk'K^2}{\pi^{3/2}},\quad\Lambda_{3,3}(3,\tau)=-\i \frac{kk'K^2}{3\pi^{5/2}}\left[\pi^2+4(2k^2-1)K^2\right].\]
\end{proposition}

To present our next results, we define the following functions
\begin{align*}
	F_{\sn}(y,k')& :=\frac{\mi}{k}\,\frac{1}{\cosh(\pi y)}\left\{k^2\int_0^{2Ky}\nd(u,k')^2\md u+2(E-K)y\right\},\\[.5ex]
	F_{\cn}(y,k')& :=\frac{1}{k}\,\frac{1}{\sinh(\pi y)}\left\{k^2\int_0^{2Ky}\nd(u,k')^2\md u+2(E-K)y\right\},\\[.5ex]
	F_{\dn}(y,k')& :=\frac{1}{\sinh(\pi y)}\left\{\int_0^{2Ky}\dn(u,k')^2\md u+2(E-K)y\right\},
\end{align*}
which are precisely the Fourier transforms $\mathcal{F}[f]$ in identities 19-21 in Table \ref{table:2}.

\begin{proposition}\label{Prop_7}
For $\im(\tau)>0$ and $\re(s)>0$ define
	\begin{align*}
		\Lambda_{0,1}(s,\tau)&:=\pi^{-\frac{s}{2}}\Gamma\Big(\frac{s}{2}\Big)\sum_{\substack{n\in\ZZ\\m\geqslant 0}}(-1)^m\frac{\big[m+\frac{1}{2}+(n+\frac{1}{2})\tau\big]^{-s}}{\sin(\pi(n+\frac{1}{2})\tau)},\\
		\Lambda_{1,0}(s,\tau)&:=\pi^{-\frac{s}{2}}\Gamma\Big(\frac{s}{2}\Big)\sum_{\substack{n\in\ZZ\sm\{0\}\\m\geqslant 0}}(-1)^{m+n}\frac{\big[m+\frac{1}{2}+n\tau\big]^{-s}}{\sin(\pi n\tau)}.
	\end{align*}
	The functions $\Lambda_{0,1}(s,\tau)$ and $\Lambda_{1,0}(s,\tau)$ can be analytically continued to entire functions of $s$ and they satisfy
	\begin{equation*}
		\Lambda_{0,1}(s,\tau)-\frac{\mi}{\tau}\,\Lambda_{1,0}\big(1-s,-1/\tau\big)=-\frac{\mi}{\tau}\frac{2}{\sqrt{\pi}}\Lambda(s-1,\chi_4), \;\;\; s\in \c.
	\end{equation*}
	Moreover, for all $n\in\ZZ_{\geqslant 0}$,
	\begin{align*}
		\Lambda_{0,1}(2n+2,\tau)& =\frac{\mi kK}{2\sqrt{\pi}}\frac{(-4 \pi)^{-n}}{\Gamma(n+\frac{3}{2})}\left[\frac{\md^{2n+1}}{\md y^{2n+1}}F_{\sn}(y,k')\right]_{y=0}, \\[.5ex]
		\Lambda_{1,0}(2n+2,\tau)& =\frac{\left|E_{2n+2}\right|}{(2n+2)!}\,\pi^{n+1}(n+1)!-\frac{k'K}{2\sqrt{\pi}}\frac{(-4 \pi)^{-n}}{\Gamma(n+\frac{3}{2})}\left[\frac{\md^{2n+1}}{\md y^{2n+1}}\frac{\sn(2Ky,k')}{\cosh(\pi y)}\right]_{y=0},
	\end{align*}
	where $\{E_{n}\}_{n\ge 0}$ are Euler numbers. In particular,
	\[\Lambda_{0,1}(2,\tau)=\frac{2K}{\pi}\left(k'^2 K-E\right),\quad\Lambda_{1,0}(2,\tau)=\frac{\pi}{2}-\frac{2}{\pi}k'K^2.\]
\end{proposition}

\begin{proposition}\label{Prop_8}
For $\im(\tau)>0$ and $\re(s)>0$ define
	\begin{equation*}
		\Lambda_{0,2}(s,\tau):=\pi^{-\frac{s}{2}}\Gamma\Big(\frac{1+s}{2}\Big)\sum_{\substack{n\in\ZZ\\m\geqslant 0}}(-1)^m\frac{\big[m+\frac{1}{2}+n\tau\big]^{-s}}{\cos(\pi n\tau)}
	\end{equation*}
	and for $\re(s)>1$ define
	\begin{equation*}
		\Lambda_{2,0}(s,\tau):=\pi^{-\frac{s}{2}}\Gamma\Big(\frac{1+s}{2}\Big)\sum_{\substack{n\in\ZZ\sm\{0\}\\m\geqslant 0}}\frac{\big[m+\frac{1}{2}+n\tau\big]^{-s}}{\sin(\pi n\tau)}.
	\end{equation*}
		The functions $\Lambda_{0,2}(s,\tau)$ and $\Lambda_{2,0}(s,\tau)$ can be analytically continued to entire functions of $s$ and they satisfy
	\begin{equation*}
		\Lambda_{0,2}(s,\tau)+\frac{\mi}{\tau}\,\Lambda_{2,0}\big(1-s,-1/\tau\big)=\frac{\mi}{\tau}\frac{s-1}{\sqrt{\pi}}\left(1-2^{2-s}\right)\Lambda(s-1), \;\;\; s\in \c.
	\end{equation*}
	Moreover, for all $n\in\ZZ_{\geqslant 0}$,
	\begin{align*}
		\Lambda_{0,2}(2n+1,\tau)& =\frac{(-4 \pi)^{-n}K}{\Gamma(n+\frac{1}{2})}\left[\frac{\md^{2n}}{\md y^{2n}}F_{\dn}(y,k')\right]_{y=0}, \\[.5ex]
		\Lambda_{2,0}(2n+1,\tau)& =\frac{(2n+1) n!}{\pi^{n+\frac{3}{2}}}\zeta\big(2n+2,\tfrac{1}{2}\big)-\frac{(-4 \pi)^{-n}K}{\Gamma(n+\frac{1}{2})}\left[\frac{\md^{2n}}{\md y^{2n}}\frac{\dn(2Ky,k')}{\cosh(\pi y)}\right]_{y=0}.
	\end{align*}
    where $\zeta(s,a)$ is the Hurwitz zeta function. In particular,
	\[\Lambda_{0,2}(1,\tau)=\frac{2KE}{\pi^{3/2}},\quad\Lambda_{2,0}(1,\tau)=\frac{\pi-2K}{2\sqrt{\pi}}.\]
\end{proposition}

\begin{proposition}\label{Prop_9}
For $\im(\tau)>0$ and $\re(s)>0$ define
	\begin{equation*}
		\Lambda_{0,3}(s,\tau):=\pi^{-\frac{s}{2}}\Gamma\Big(\frac{1+s}{2}\Big)\sum_{\substack{n\in\ZZ\\m\geqslant 0}}(-1)^m\frac{\big[m+\frac{1}{2}+(n+\frac{1}{2})\tau\big]^{-s}}{\cos(\pi(n+\frac{1}{2})\tau)}
	\end{equation*}
	and for $\re(s)>1$ define
	\begin{equation*}
		\Lambda_{3,0}(s,\tau):=\pi^{-\frac{s}{2}}\Gamma\Big(\frac{1+s}{2}\Big)\sum_{\substack{n\in\ZZ\sm\{0\}\\m\geqslant 0}}(-1)^n\frac{\big[m+\frac{1}{2}+n\tau\big]^{-s}}{\sin(\pi n\tau)}.
	\end{equation*}
	The functions $\Lambda_{0,3}(s,\tau)$ and $\Lambda_{3,0}(s,\tau)$ can be analytically continued to entire functions of $s$ and they satisfy
	\begin{equation*}
		\Lambda_{0,3}(s,\tau)+\frac{\mi}{\tau}\,\Lambda_{3,0}\big(1-s,-1/\tau\big)=\frac{\mi}{\tau}\frac{s-1}{\sqrt{\pi}}\left(1-2^{2-s}\right)\Lambda(s-1), \;\;\; s\in \c. 
	\end{equation*}
	Moreover, for all $n\in\ZZ_{\geqslant 0}$,
	\begin{align*}
		\Lambda_{0,3}(2n+1,\tau)& =kK\frac{(-4 \pi)^{-n}}{\Gamma(n+\frac{1}{2})}\left[\frac{\md^{2n}}{\md y^{2n}}F_{\cn}(y,k')\right]_{y=0},\\[.5ex]
		\Lambda_{3,0}(2n+1,\tau)& =\frac{(2n+1) n!}{\pi^{n+\frac{3}{2}}}\zeta\big(2n+2,\tfrac{1}{2}\big)-k'K\frac{(-4 \pi)^{-n}}{\Gamma(n+\frac{1}{2})}\left[\frac{\md^{2n}}{\md y^{2n}}\frac{\cn(2Ky,k')}{\cosh(\pi y)}\right]_{y=0}. 
	\end{align*}
	In particular,
	\[\Lambda_{0,3}(1,\tau)=\frac{2K}{\pi^{3/2}}\left(E-k'^2 K\right),\quad\Lambda_{3,0}(1,\tau)=\frac{\pi-2k'K}{2\sqrt{\pi}}.\]
\end{proposition}

\begin{proposition}\label{Prop_10}
    For $\re(s)>0$ define
	\begin{equation*}
		\Lambda_{0,0}(s,\tau):=\pi^{-\frac{s}{2}}\Gamma\Big(\frac{s}{2}\Big)\sum_{\substack{n\in\ZZ\sm\{0\}\\m\geqslant 0}}(-1)^m\frac{\big[m+\frac{1}{2}+n\tau\big]^{-s}}{\sin(\pi n\tau)}.
	\end{equation*}
	The function $\Lambda_{0,0}(s,\tau)$ can be analytically continued to an entire function of $s$ and it satisfies
	\begin{equation*}
		\Lambda_{0,0}(s,\tau)-\frac{\mi}{\tau}\,\Lambda_{0,0}\big(1-s,-1/\tau\big)=\frac{2}{\sqrt{\pi}}\Big(\Lambda(s+1,\chi_4)-\frac{\mi}{\tau}\,\Lambda(s-1,\chi_4)\Big), \;\;\;
		s\in \c.
	\end{equation*}
	Moreover, for all $n\in\ZZ_{\geqslant 0}$,
	\begin{equation*}
		\Lambda_{0,0}(2n+2,\tau)=-\frac{\sqrt{\pi}}{4K'}\frac{(-4 \pi)^{-n}}{\Gamma(n+\frac{3}{2})}\left[\frac{\md^{2n+1}}{\md y^{2n+1}}\frac{2Ky+K\,\mathfrak{f}(y,k')-K'\tanh(\pi y)}{\cosh(\pi y)}\right]_{y=0},
	\end{equation*}
	where ${\mathfrak f}$ is defined in 
	\eqref{def_mathfrak_f}. 
	In particular,
	$$
	\Lambda_{0,0}(2,\tau)=\frac{1}{2\pi}(\pi^2-4KE). 
	$$
\end{proposition}

The proofs of these ten propositions are based on the following two theorems. The methods underlying these theorems are not new -- they were used, for example, in \cite{TiHB} to establish analyticity properties and the functional equation for the Riemann zeta function. We denote the Mellin transform of a function $f:[0,\infty)\to\CC$ by
\[\mathcal{M}[f](z):=\int_0^{\infty}f(x)x^{z-1}\md x\]
and the Fourier cosine and sine transforms by
\[\mathcal{F}_c[f](y):=\int_0^{\infty}f(x)\cos(xy)\md x,\quad\mathcal{F}_s[f](y):=\int_0^{\infty}f(x)\sin(xy)\md x.\]

\begin{theorem}\label{thm: cosine transform}
	Assume that $f:\RR\to\CC$ is a smooth even (respectively, odd) periodic function, and let $g(x):=f(x)/\cosh(\pi x)$ (respectively, $g(x):=f(x)/\sinh(\pi x)$). Denote $g_c(y):=\mathcal{F}_c[g](y)$. The Mellin transforms $G(z):=\mathcal{M}[g](z)$ and $G_c(z):=\mathcal{M}\left[g_c\right](z)$ can be analytically continued to meromorphic functions that have simple poles at nonpositive even integers and satisfy the following properties:
	\begin{itemize}
			\item[{\rm(i)}] $\Gamma(z)\cos\left(\frac{\pi}{2}z\right)G(1-z)=G_c(z)$ for $z\in\CC$;
		\item[{\rm(ii)}] for $n \in {\mathbb Z}_{\geqslant 0}$ 
		\begin{align*}
		G(2n+1)&= (-1)^n (2n)! \times \Res\limits_{z=-2n }G_c(z)= (-1)^n g_c^{(2n)}(0), \\
		\frac{2}{\pi} G_c(2n+1)&=  (-1)^n (2n)! \times \Res\limits_{z=-2n}G(z)=  (-1)^n g^{(2n)}(0).				
		\end{align*}
	\end{itemize}
\end{theorem}

\begin{theorem}\label{thm: sine transform}
	Let $f:\RR\to\CC$ be a smooth even (respectively, odd) periodic function that satisfies $f(0)=0$, and let $g(x):=f(x)/\sinh(\pi x)$ (respectively, $g(x):=f(x)/\cosh(\pi x)$). Denote $g_s(y):=\mathcal{F}_s[g](y)$. The Mellin transforms $G(z):=\mathcal{M}[g](z)$ and $G_s(z):=\mathcal{M}\left[g_s\right](z)$ can be analytically continued to meromorphic functions that have simple poles at negative odd integers and satisfy the following properties:
	\begin{itemize}
			\item[{\rm(i)}] $\Gamma(z)\sin\left(\frac{\pi}{2}z\right)G(1-z)=G_s(z)$ for $z\in\CC$;
		\item[{\rm(ii)}] for $n \in {\mathbb Z}_{\geqslant 0}$
	\begin{align*}	
		G(2n+2)&= (-1)^n (2n+1)! \times \Res\limits_{z=-2n-1}G_s(z)= (-1)^n g_s^{(2n+1)}(0),\\
				\frac{2}{\pi}  G_s(2n+2)&=  (-1)^n (2n+1)! \times \Res\limits_{z=-2n-1}G(z)=  (-1)^n g^{(2n+1)}(0).
	\end{align*}			 
	\end{itemize}
\end{theorem}

We will only prove Theorem \ref{thm: cosine transform}, as the proof of Theorem \ref{thm: sine transform} follows exactly the same steps.

\vspace{0.25cm}

\noindent
{\bf Proof of Theorem \ref{thm: cosine transform}}:
Consider the case when $f$ is odd. Since we assumed that $f$ is odd, smooth, and periodic, it follows that $f(0)=0$ and $g(x)=f(x)/\sinh(x)$ is an even Schwartz-class function. The Fourier transform of $g$ is also a Schwartz-class function, and since $g$ is even, its Fourier transform is equal (up to a constant) to its cosine transform $g_c(y)$. Thus, $g_c(y)$ is an even Schwartz-class function. The fact that $G(2n+1)=(-1)^n g_c^{(2n)}(0)$ follows by differentiating both sides of the integral identity
\[
g_c(y)=\int_0^{\infty}g(x)\cos(xy)\d x
\]
with respect to $y$ and setting $y=0$.

Next, we write
\[
G(z):=\mathcal{M}[g](z)=\int_0^1 g(x)x^{z-1}\md x+\int_1^{\infty}g(x)x^{z-1}\md x=G_1(z)+G_2(z).
\]
Denoting $a_j:=g^{(2j)}(0)/(2j)!$, we have for every $n \in {\mathbb N}$,
\[
G_1(z)=\sum_{j=0}^n\frac{a_j}{z+2j}+\int_0^1\Big[g(x)-\sum_{j=0}^n a_j x^{2j}\Big]x^{z-1}\md x,
\]
and the integral on the right-hand side of the above equation is an analytic function of $z$ in the half-plane $\re(z)>-2n-2$. Thus, $G_1$ is a meromorphic function with simple poles at $z=-2j$ and residues $\Res_{z=-2j}G_1(z)=a_j$. The same applies to $G(z)=G_1(z)+G_2(z)$, since $G_2$ is entire (which follows from the fact that $g(x)$ decays exponentially as $x\to +\infty$).

Therefore, we have proved that $G(z)$ can be continued to a meromorphic function having simple poles at $z=-2n$ for $n \in {\mathbb Z}_{\geqslant 0}$. The same applies to $G_c(z)$ (since, as we noted above, $g_c$ is an even Schwartz function). The functional equation $\Gamma(z) \cos\left(\frac{\pi}{2}z\right) G(1-z)= G_c(z)$ now follows from relation (b$'$) on page 3 in \cite{Ober1}.

We established above that 
\begin{equation}\label{residue_value_G}
\Res_{z=-2n}G(z)=g^{(2n)}(0)/(2n)!,  \;\;\; G(2n+1)=(-1)^n g_c^{(2n)}(0),
\end{equation}
for $n \in {\mathbb Z}_{\geqslant 0}$. Using the functional equation established in item (i), we derive 
\[
\Res_{z=-2n}G_c(z)=g_c^{(2n)}(0)/(2n)!,  \;\;\; G_c(2n+1)=\frac{\pi}{2}(-1)^n g^{(2n)}(0).
\]
One could also obtain the above relations from \eqref{residue_value_G} by noting that ${\mathcal F}_c [g_c](y)=\frac{\pi}{2} g(y)$.

This completes the proof of Theorem \ref{thm: cosine transform} in the case when $f$ is odd. The proof for the case when $f$ is even follows the same steps and is left to the reader.
\qed

Theorems \ref{thm: cosine transform} and \ref{thm: sine transform} show that Fourier cosine and sine transform pairs lead (via the Mellin transform) to meromorphic functions satisfying a functional equation -- a reflection formula under the transformation $s\mapsto 1-s$. The next table shows which of the $\Lambda_{j,l}(s,\tau)$ functions are connected with which Fourier transform pairs.

\begin{table}[htbp]
	\centering
	\begin{tabular}{|c|c||c|c||c|c|}
		\hline
		\rule{0pt}{1\normalbaselineskip}
		\text{Double series} & \text{Identity} & \text{Double series} & \text{Identity} & \text{Double series} & \text{Identity}\\[.5ex]
		\hline
		\rule{0pt}{1\normalbaselineskip}
		        $\Lambda_{0,0}$ & 25 &  &  &  & \\[1ex]
		         $\Lambda_{1,1}$ & 4 & $\Lambda_{2,2}$& 17 & $\Lambda_{3,3}$ & 15\\[1ex]       
        $\Lambda_{0,1},\Lambda_{1,0}$ & 21 & $\Lambda_{0,2},\Lambda_{2,0}$& 20 & $\Lambda_{0,3},\Lambda_{3,0}$ & 19 \\[1ex]
		$\Lambda_{1,2}, \Lambda_{2,1}$ & 5, 18 & $\Lambda_{1,3}, \Lambda_{3,1}$ & 6, 13 & $\Lambda_{2,3}, \Lambda_{3,2}$ & 14, 16\\
		\hline
	\end{tabular}
		\caption{Eisenstein-type series $\Lambda_{j,l}(s,\tau)$  and the corresponding Fourier transform identity. }\label{table:3}
\end{table}

Our plan is to illustrate the use of Theorems \ref{thm: cosine transform} and \ref{thm: sine transform} by proving Proposition \ref{Prop_3}.
First, we recall some facts about the Hurwitz zeta function and its alternating version (see \cite[\S 25.11]{NIST} and \cite{HU2024,Nemes_2017}), defined as 
$$
\zeta(s,a):=\sum\limits_{n\geqslant 0} (n+a)^{-s}, \;\;\; 
\zeta_E(s,a):=\sum\limits_{n\geqslant 0} (-1)^n (n+a)^{-s}. 
$$
Here $|\arg(a)|<\pi$, and we require $\re(s)>1$ ($\re(s)>0$) for the first (respectively, second) infinite series to converge.
For $\re(a)>0$, the following integral representations hold:
\begin{align}
\label{Hurwitz_integral}
\zeta(s,a)=\frac{1}{\Gamma(s)}
\int_0^{\infty} \frac{x^{s-1} e^{-ax}}{1-e^{-x}} \d x, \;\;\; \re(s)>1,  \\
\label{Hurwitz2_integral}
\zeta_E(s,a)=\frac{1}{\Gamma(s)}
\int_0^{\infty} \frac{x^{s-1} e^{-ax}}{1+e^{-x}} \d x, \;\;\; \re(s)>0. 
\end{align}
It is known that $\zeta(s,a)$ can be extended to an analytic function of $(a,s)$ in the domain $|\arg(a)|<\pi$, $s\in \c \setminus \{1\}$, and as a function of $s$, it has a simple pole at $s=1$ with residue $1$. The alternating Hurwitz zeta function can be extended to an analytic function of $(a,s)$ in the domain 
$|\arg(a)|<\pi$, $s\in \c$. 

The asymptotics of $\zeta(s,a)$ and $\zeta_E(s,a)$ for large $a$ were studied in \cite{HU2024, Nemes_2017}. In particular, from \cite[Theorem 1.2]{Nemes_2017} it follows that for any fixed $\epsilon \in (0,\pi)$, as $a\to \infty$ in the sector $|\arg(a)|<\pi-\epsilon$
we have	
\begin{equation}\label{zeta_asymptotics}
\zeta(s,a)=\frac{1}{2} a^{-s} + \frac{a^{1-s}}{s-1} + O(|a|^{-s-1}),  
\end{equation}
and this holds uniformly in $s$ on compact subsets of $\{\re(s)>-1\} \setminus \{1\}$. The corresponding result for 
the alternating Hurwitz zeta function (see \cite[Theorem 3.7]{HU2024}) states that for any fixed $\epsilon \in (0,\pi)$, as $a\to \infty$ in the sector $|\arg(a)|<\pi-\epsilon$
we have	
\begin{equation}\label{zeta_E_asymptotics}
\zeta_E(s,a)=\frac{1}{2}a^{-s} 
+ \frac{s}{4} a^{-s-1} + O(|a|^{-s-3}), 
\end{equation}
uniformly in $s$ on compact subsets of $\{\re(s)>-3\}$.

\vspace{0.25cm}
\noindent
{\bf Proof of Proposition \ref{Prop_3}:} 
We assume first that $k \in (0,1)$, so that $\tau=\i K'/K$ satisfies $\re(\tau)=0$ and $\im(\tau)>0$.
We rewrite definitions \eqref{def_Lambda13} and \eqref{def_Lambda31} of $\Lambda_{1,3}$ and $\Lambda_{3,1}$ in the form
	\begin{equation}\label{def_Lambda_13_new}
		\Lambda_{1,3}(s,\tau)=\pi^{-\frac{s}{2}}\Gamma\Big(\frac{s}{2}\Big)\sum_{n\in\ZZ}\left(-1\right)^{n}\frac{\zeta_E\big(s,\frac{1}{2} + (n+\frac{1}{2})\tau\big)}{\cos(\pi (n+\frac{1}{2})\tau)}, \;\;\; \re(s)>0, 
	\end{equation}
	and 
	\begin{equation}\label{def_Lambda_31_new}
		\Lambda_{3,1}(s,\tau):=\pi^{-\frac{s}{2}}\Gamma\Big(\frac{s}{2}\Big)\sum_{n\in\ZZ}(-1)^n\frac{\zeta\big(s,\frac{1}{2} + (n+\frac{1}{2})\tau\big)}{\sin(\pi (n+\frac{1}{2} )\tau )}, \;\;\; \re(s)> 1. 
	\end{equation}

Denote $f(x)=\cd(2K'x,k)-1$. We verify that $f(x)$ is a smooth, even, periodic function satisfying $f(0)=0$. As in Theorem \ref{thm: sine transform}, we define $g(x)=f(x)/\sinh(\pi x)$. 
Fourier series in \cite[Eq. (22.11.4)]{NIST} gives
	\begin{align*}
		\cd(2K'x,k) =\frac{\pi\mi}{Kk}\sum_{n\geqslant 0}(-1)^n\frac{\cos((2n+1)\pi\frac{K'}{K}x)}{\sin(\pi(n+\frac{1}{2})\tau)}, \;\;\; x\in \r. 
	\end{align*}
Thus, for $\re(z)>1$ we have
	\begin{align}
		G(z)=\mathcal{M}[g](z)&=\frac{\pi\mi}{Kk}\sum_{n=0}^{\infty}\frac{(-1)^n}{\sin(\pi(n+\frac{1}{2})\tau)}\int_0^{\infty}\frac{\cos((2n+1)\pi\frac{K'}{K}x)}{\sinh(\pi x)}x^{z-1}\md x
		- \int_0^{\infty} \frac{x^{z-1}\d x}{\sinh(\pi x)} 		
		 \nonumber\\
		& =\frac{\pi\mi}{Kk}\frac{\Gamma(z)}{(2\pi)^z}\sum_{n\geqslant 0}\frac{(-1)^n}{\sin(\pi(n+\frac{1}{2})\tau)} 
		\Big[\zeta(z,\tfrac{1}{2}+(n+\tfrac{1}{2})\tau)+
		\zeta(z,\tfrac{1}{2}-(n+\tfrac{1}{2})\tau)\Big]
		\label{relation_G_Lambda_3,1}\\
		&\qquad \qquad \qquad \qquad \qquad \qquad -2\left(1-2^{-z}\right)\pi^{-z}\Gamma(z)\zeta(z) 
		 \nonumber\\ \nonumber
		& =\frac{\mi}{2Kk}\,\pi^{\frac{1-z}{2}}\Gamma\Big(\frac{1+z}{2}\Big)\Lambda_{3,1}(z,\tau)-2\left(1-2^{-z}\right)\pi^{-z}\Gamma(z)\zeta(z). 
	\end{align}
Here we used Lebesgue’s dominated convergence theorem to interchange the summation and integration. In deriving the above formula, we also used \eqref{sinh_integral_zeta},  \eqref{Hurwitz_integral}, and \eqref{def_Lambda_31_new}.  Theorem \ref{thm: sine transform} tells us that $G(z)$ extends to a meromorphic function on $\c$, whose only poles are simple and lie at negative odd integers. The Riemann zeta function $\zeta(z)$ is analytic on $\c$, except for a simple pole at $z=1$, and $\Gamma(z)$ has simple poles at $z \in {\mathbb Z}_{\leqslant 0}$. These facts, together with formula \eqref{relation_G_Lambda_3,1}, imply that $\Lambda_{3,1}(z,\tau)$ extends to a meromorphic function on $\CC$, which has a simple pole at $z=1$ with residue $-2\mi kK/\pi$. 
	
Fourier transform pair 13 in Table \ref{table:1} and identity \eqref{Fourier/sinh} imply 
	\begin{equation}\label{g(x) sine transform}
		g_s(y)=\int_0^{\infty}g(x)\sin(xy)\md x=- \frac{k'}{2}\frac{ \sd(Ky/\pi,k')}{\cosh(y/2)}, \;\;\; y\in \r. 
	\end{equation}	
In the same manner as above for $\mathcal{M}[g](z)$, we compute the Mellin transform of $g_s(y)$ and obtain
\begin{equation}	\label{G_s_prop3}
	G_s(z)=\mathcal{M}\left[g_s \right](z)=
	-\frac{\i}{K' k} 2^{z-2} \pi^{\frac{1+z}{2}} 
	\Gamma\Big(\frac{1+z}{2}\Big)\Lambda_{1,3}(z,-1/\tau), 
	\;\;\; \re(z)>0.
\end{equation}	
Expressions \eqref{relation_G_Lambda_3,1} and \eqref{G_s_prop3}, together with the functional equation in Theorem \ref{thm: sine transform}(i), after some simplification, lead to the functional equation \eqref{reflection formula for Lambda_1,3 and Lambda_3,1}. Formulas \eqref{Lambda_1,3 values} and \eqref{Lambda_3,1 values}, which give the values of $\Lambda_{1,3}(2n+2,\tau)$ and $\Lambda_{3,1}(2n+2,\tau)$, follow from Theorem \ref{thm: sine transform}(ii). The special cases $n=0,1$ can be computed with the help of \cite[\S 22.5(i) and \S 22.13]{NIST}. 
	
So far, we have established Proposition \ref{Prop_3} in the special case when $\re(\tau)=0$. We now extend these results to $\im(\tau)>0$. 
In light of \eqref{zeta_E_asymptotics}, it is clear that the series
\eqref{def_Lambda_13_new} for $\Lambda_{1,3}$ converges uniformly in $\tau$ and $s$ on compact subsets of $\im(\tau)>0$ and $\re(s)>-3$, thus 
$\Lambda_{1,3}(s,\tau)$ is an analytic function of $(s,\tau)$ in this domain. Similarly, we establish that $\Lambda_{3,1}(s,\tau)$ is an analytic function of $(s,\tau)$ in the domain $\im(\tau)>0$ and $\{\re(s)>-1, s\neq 1\}$. Hence, by analytic continuation in $\tau$, we conclude that all statements in Proposition \ref{Prop_3} hold for $\im(\tau)>0$. \qed

The remaining Propositions 1-2 and 4-10 are established in exactly the same way. Perhaps the only nontrivial step is in establishing the formula for $\Lambda_{0,0}(2,\tau)$ in 
Proposition \ref{Prop_10}. This requires an expression for ${\mathfrak f}_y(0,k')$, which can be obtained from the identity
$$
\frac{\theta_4''(0|\tau)}{\theta_4(0|\tau)}=
8 \sum\limits_{n\geqslant 1} \frac{q^{2n-1}}{(1-q^{2n-1})^2}=\frac{4}{\pi^2}K (K-E),
$$
see the formula for ${\textnormal{III}}_2(c)$ in \cite[Eq. (30)]{LING1988}.

\section{Concluding remarks}\label{section_conclusion}

We would like to conclude by discussing two problems that we did not resolve in the present paper. First, we computed only the Fourier transform of the logarithmic derivative of $\theta_4(\cdot|\tau)$ (our identity 25). It remains to compute the Fourier transforms of the logarithmic derivatives of the theta functions $\theta_j(\cdot|\tau)$ for $j=1,2,3$. 

Secondly, a look at Table \ref{table:3} reveals that we used only thirteen Fourier transform pair identities (out of a total of 25) to derive our results concerning the sixteen Eisenstein-type series $\zeta_{j,l}(s,\tau)$ defined in \eqref{def_zeta_general2}. The following twelve Fourier transform pairs were not used: 1-3, 7-9, 10-12, and 22-24. From Tables \ref{table:1} and \ref{table:2}, we see that these identities correspond to Fourier transforms of the six Jacobi functions $*$c$(2K'x,k)$ and $*$s$(2K'x,k)$ (where $*$ stands for the letter c, d, n, or s). These are precisely the Jacobi functions that have poles on the lattice $(m+\frac{1}{2})K/K'+n \i$ or $mK/K'+n \i$, with $m,n \in {\mathbb Z}$ (see \cite[Table 22.4.1]{NIST}). 

Thus, in these twelve cases (Fourier transform pairs 1-3, 7-9, 10-12, and 22-24), the function appearing as the integrand in the Fourier transform has infinitely many poles on the real line. These poles prevent us from directly applying Theorems \ref{thm: cosine transform} and \ref{thm: sine transform}. However, we believe that by a suitable deformation of the contour of integration in the definitions of the Mellin and Fourier cosine/sine transforms, one could derive modified versions of Theorems \ref{thm: cosine transform} and \ref{thm: sine transform} that would apply to Fourier transform identities 1-3, 7-9, 10-12, and 22-24. 

This should lead to analogues of our results in Section \ref{section_Mellin} for a new family of Eisenstein-type series of the form 
 \begin{equation}\label{extended_zeta}
	\sideset{}{'}\sum_{\substack{n\in\ZZ\\m\geqslant 1}}\epsilon_1^m \epsilon_2^n \frac{\big [m+(n+\frac{d}{2})\tau \big]^{-s}}{\sin(\pi(\frac{c}{2}+(n+\frac{d}{2})\tau))},
\end{equation}
where $\epsilon_i \in \{-1,1\}$ and $c,d \in \{0,1\}$. Further evidence that such an extension of our results is possible comes from the work of Tsumura \cite{Tsumura2008,TSUMURA_2009,Tsumura2010}, who  evaluated in closed form certain special cases (with $s \in {\mathbb N}$ and $\tau=\i$) of the double series in \eqref{extended_zeta}. We leave the investigation of this new family of Eisenstein-type series to future work.

\section*{Acknowledgements}
The research was supported by the Natural Sciences and Engineering Research Council of Canada. We are grateful to an anonymous referee for carefully reading the paper and for providing very helpful comments.


\end{document}